\documentclass{elsart3-1}

\usepackage{amsmath,amssymb,graphicx,epsfig,color}
\usepackage[english,francais]{babel}

\newtheorem{theorem}{Theorem}[section]

\newtheorem{e-proposition}[theorem]{Proposition}

\newtheorem{e-definition}[theorem]{Definition\rm}

\renewcommand{\theequation}{\arabic{equation}}
\def\zz{\mathbb{Z}}
\def\rr{\mathbb{R}}
\def\nn{\mathbb{N}}

\def\og{\leavevmode\raise.3ex\hbox{$\scriptscriptstyle\langle\!\langle$~}}
\def\fg{\leavevmode\raise.3ex\hbox{~$\!\scriptscriptstyle\,\rangle\!\rangle$}}

\begin{document}

\begin{frontmatter}

\selectlanguage{english}
\title{Localized solutions for the finite difference semi-discretization of the wave equation}

\vspace{-2.6cm}

\selectlanguage{francais}
\title{Solutions localis\'{e}es pour la semi-discr\'{e}tisation par diff\'{e}rences finies de l'\'{e}quation des ondes}

\selectlanguage{english}
\author[BCAM]{Aurora Marica},
\ead{marica@bcamath.org}
\author[Iker,BCAM]{Enrique Zuazua}
\ead{zuazua@bcamath.org}

\address[Iker]{Ikerbasque, Basque Foundation for Science, Alameda Urquijo 36-5, Plaza Bizkaia, 48011, Bilbao, Basque Country, Spain}

\address[BCAM]{BCAM - Basque Center for Applied Mathematics, Bizkaia Technology Park 500, 48160, Derio, Basque Country, Spain}

\begin{abstract}
We study the propagation properties of the solutions of the finite-difference space semi-discrete wave equation on an uniform grid of the whole Euclidean space. We provide a construction of high frequency wave packets that propagate along the corresponding bi-characteristic rays of Geometric Optics with a group velocity arbitrarily close to zero. Our analysis is motivated by control theoretical issues. In particular, the continuous wave equation has the so-called observability property: for a sufficiently large time, the total energy of its solutions can be estimated in terms of the energy concentrated in the exterior of a compact set. This fails to be true, uniformly on the mesh-size parameter, for the semi-discrete schemes and the observability constant blows-up at an arbitrarily large polynomial order. Our contribution consists in providing a rigorous derivation of those wave packets and in analyzing their behavior near that ray, by taking into account the subtle added dispersive effects that the numerical scheme introduces.

\vskip 0.5\baselineskip

\selectlanguage{francais}
\noindent{\bf R\'esum\'e
 } \vskip 0.5\baselineskip \noindent

On \'{e}tudie les propri\'{e}t\'{e}s de propagation des solutions de l'\'{e}quation des ondes semi-discretis\'{e}e en espace par diff\'{e}rences finies sur une grille uniforme dans tout l'espace euclidien. On r\'{e}alise une construction de paquets d'ondes concentr\'{e}s \`{a} haute fr\'{e}quence qui se propagent le long des rayons bicaract\'{e}ristiques de l'Optique G\'{e}om\'{e}trique \`{a} une vitesse de groupe arbitrairement petite. Notre analyse est motiv\'{e}e par la th\'{e}orie du contr\^{o}le. Plus pr\'{e}cisement, l'\'{e}quation des ondes continue v\'{e}rifie la propri\'{e}t\'{e} d'observabilit\'{e}: pour un temps suffisament grand, l'\'{e}nergie totale de ses solutions peut \^{e}tre estim\'{e}e  en fonction  de l'\'{e}nergie localis\'{e}e \`{a} l'ext\'{e}rieur d'un ensemble compact. Cette propri\'{e}t\'{e} n'est pas verifi\'{e}e de mani\`{e}re uniforme par rapport au pas de discr\'{e}tisation pour le sch\'{e}ma semi-discret pour un temps fini quelconque, si bien que la constante d'observabilit\'{e} semi-discr\`{e}te diverge avec un taux polynomial arbitraire. Nous donnons une construction pr\'{e}cise de ces paquets d'ondes et decrivons l'effet dispersif rajout\'{e} que le sch\'{e}ma num\'{e}rique g\'{e}n\`{e}re autour du rayon de propagation.
\end{abstract}
\end{frontmatter}

\selectlanguage{francais}
\section*{Version fran\c{c}aise abr\'eg\'ee}

Il est bien connu (cf. \cite{ZuaUnbDom}) que les solutions du probl\`{e}me de Cauchy pour l'\'{e}quation des ondes continue $d$-dimensionnelle v\'{e}rifient la propri\'{e}t\'{e} d'observabilit\'{e}: pour un temps suffisament grand, l'\'{e}nergie totale des solutions peut \^{e}tre estim\'{e}e en fonction de l'\'{e}nergie localis\'{e}e \`{a} l'ext\'{e}rieur d'un ensemble compact. Par la m\'{e}thode d'unicit\'{e} de Hilbert (HUM) introduit dans \cite{LioI}, cette propri\'{e}t\'{e} est equivalente \`{a} un r\'{e}sultat de contr\^{o}labilit\'{e} exacte par un contr\^{o}le localis\'{e} dans le compl\'{e}mentaire du m\^{e}me ensemble compact.

Dans cet article, on consid\`{e}re le probl\`{e}me de Cauchy associ\'{e} \`{a} l'\'{e}quation des ondes $d$-dimensionnelle semi-discretis\'{e}e en espace par un sch\'{e}ma centr\'{e} en diff\'{e}rences finies dans un maillage uniforme. Nous avons deux objectifs. Le premier est de montrer l'\'{e}xistence de paquets d'ondes concentr\'{e}s autour d'un nombre d'onde fix\'{e} a priori, \`{a} haute fr\'{e}quence, propageant avec une vitesse de groupe arbitrairement petite autour des rayons bicharact\'{e}ristiques du sch\'{e}ma semi-discret et donc peu visibles depuis le domaine d'observation. On d\'{e}duit dans \cite{Erv} que la constante d'observabilit\'{e} semi-discr\`{e}te diverge au moins d'une mani\`{e}re polynomiale arbitraire, et ce pour tout temps d'observabilit\'{e} fini. Le deuxi\`{e}me objectif est de donner une forme asymptotique pr\'{e}cise pour ces paquets d'ondes, de fa\c{c}on \`{a} mettre en \'{e}vidence la dispersion num\'{e}rique qui n'appara\^{i}t pas dans le mod\`{e}le continu.

Notre r\'{e}sultat compl\`{e}te la litt\'{e}rature existante sur l'observabilit\'{e} et le contr\^{o}le des sch\'{e}mas d'approximation de l'\'{e}quation des ondes. En particulier, dans \cite{ZuaPOC}, pour la semi-discr\'{e}tisation de l'\'{e}quation des ondes par diff\'{e}rences finies dans des domaines born\'{e}s comme les cubes $d$-dimensionnels, il a \'{e}t\'{e} prouv\'{e} que l'in\'{e}galit\'{e} d'observabilit\'{e} explose au sens o\`{u} la constante d'observabilit\'{e} explose lorsque le pas du maillage tend vers z\'{e}ro. Dans  \cite{MicuBio}, on montre aussi que pour l'\'{e}quation des ondes unidimensionnelle sur un intervalle born\'{e} observ\'{e}e par le bord, la constante d'observabilit\'{e} diverge de mani\`{e}re exponentielle. Ce type de ph\'{e}nom\`{e}ne pathologique appara\^{i}t aussi au niveau des estimations dispersives de Strichartz associ\'{e}es \`{a} l'\'{e}quation de Schr\"{o}dinger semi-discretis\'{e}e en espace par diff\'{e}rences finies, cf. \cite{LivEZDispSchr}.

L'avantage de notre construction de paquets d'ondes num\'{e}riques \`{a} haute fr\'{e}quence est qu'elle peut \^{e}tre g\'{e}n\'{e}ralis\'{e}e \`{a} d'autres approximations plus sophistiqu\'{e}es de l'\'{e}quation des ondes comme celles obtenues par les m\'{e}thodes de Galerkin discontinus ou par des \'{e}l\'{e}ments finis classiques d'ordre sup\'{e}rieur, o\`{u} le symbole ``num\'{e}rique" en Fourier du Laplacien est une matrice donnant lieu \`{a} plusieurs relations de dispersion, cf. \cite{MarPhD}.

En outre, les paquets d'ondes que nous construisons sont aussi localis\'{e}s dans l'espace physique et peuvent donc \^{e}tre adapt\'{e}s pour des probl\`{e}mes aux limites dans des domaines born\'{e}s avec des conditions aux limites diverses.

\selectlanguage{english}
\section{Introduction and problem formulation}
For a finite time $T>0$, consider the Cauchy problem associated to the $d$-dimensional wave equation:
\begin{equation}\left\{\begin{array}{ll}\partial^2_t\phi-\Delta \phi=0,&x\in\rr^d,t\in(0,T]
\\\phi(x,0)=\phi^0(x),\partial_t\phi(x,0)=\phi^1(x),&x\in\rr^d.\end{array}\right.\label{ContWave}\end{equation}
The problem (\ref{ContWave}) is well posed in $\dot{H}^1(\rr^d)\times L^2(\rr^d)$ ($\dot{H}^1(\rr^d)$ is the completion of $C_c^{\infty}(\rr^d)$ with respect to the semi-norm $\|\cdot\|_{\dot{H}^1(\rr^d)}^2=\|\nabla\cdot\|_{L^2(\rr^d)}^2$). The total energy of its solutions is conserved in time and is defined by $E(\phi^0,\phi^1)=(\|\phi(t)\|^2_{\dot{H}^1(\rr^d)}+\|\partial_t\phi(t)\|^2_{L^2(\rr^d)})/2$.

For all finite time $T>2$, there exists a constant $C(T)>0$ such that, for all $(\phi^0,\phi^1)\in \dot{H}^1(\rr^d)\times L^2(\rr^d)$, the following observability inequality holds (see \cite{ZuaUnbDom}):
\begin{equation}E(\phi^0,\phi^1)\leqslant C(T)\int\limits_0^T\int\limits_{\Omega^d}|\partial_t\phi(x,t)|^2\,dx\,dt.\label{ContObsIneq}\end{equation}
 The observation region $\Omega^d:=\rr^d\setminus B^d(0,1)$ is the complement of the $d$-dimensional unit ball. This observability problem is motivated by controllability issues, since, by means of the Hilbert Uniqueness Method (HUM) introduced in \cite{LioI},  the inequality (\ref{ContObsIneq}) is equivalent  to the fact that, for all $T>2$ and all initial data $(u^0,u^1)\in \dot{H}^1(\rr^d)\times L^2(\rr^d)$, there exists a control function $f \in L^2(\Omega^d \times (0, T))$ such that the solution of the inhomogeneous Cauchy problem
\begin{equation}\left\{\begin{array}{ll}\partial^2_tu-\Delta u=f\chi_{\Omega^d}(x),&x\in\rr^d,\ t\in(0,T]\\u(x,0)=u^0(x),\ \partial_tu(x,0)=u^1(x),& x\in\rr^d\end{array}\right.\label{controlled}\end{equation}
satisfies $u(x,T)=\partial_tu(x,T)=0$ for all $x\in\rr^d$, where $\chi_{\Omega^d}$ is the characteristic function of the set $\Omega^d$.

These issues are by now well understood for the continuous wave
equation. In particular, observability inequalities of the form (\ref{ContObsIneq}) and several variants
hold under the sharp Geometric Control Condition (GCC) (cf. \cite{BaLebRa}) requiring that all rays of Geometric Optics enter the observation subdomain during the observability time. As shown in \cite{Ral}, when the GCC is not satisfied, the observability property fails
because of the existence of Gaussian beam solutions localized around a
bi-characteristic ray that escape the observation region during the time interval $[0,T]$.

This paper is devoted to analyze these issues for the finite difference semi-discrete schemes. Given a mesh size $h>0$, we define an uniform grid of the whole Euclidean space by $x_{\mathbf{j}}=h\mathbf{j}$, $\mathbf{j}\in\zz^d$. We also introduce two discrete operators: the gradient $\nabla_h^+=(\partial_{h,k}^+)_{k=1,\cdots,d}$ and the Laplacian $\Delta_h$, with
\begin{equation}\partial_{h,k}^+\overrightarrow{f}=\frac{\overrightarrow{f}_{\cdot+\mathbf{e}_k}-\overrightarrow{f}}{h},\quad \Delta_h\overrightarrow{f}=\frac{1}{h^2}\sum\limits_{l=1}^d(\overrightarrow{f}_{\cdot+\mathbf{e}_l}-2\overrightarrow{f}+
\overrightarrow{f}_{\cdot-\mathbf{e}_l}),\label{DiscOper}\end{equation}
where $\overrightarrow{f}=(f_{\mathbf{j}})_{\mathbf{j}\in\zz^d}$ is any sequence and $(\mathbf{e}_l)_{l=1}^d$ is the canonical basis in $\rr^d$.

We consider the finite difference space semi-discretization of the wave equation (\ref{ContWave}):
\begin{equation}\left\{\begin{array}{ll}\partial^2_t\phi_{\mathbf{j}}^h(t)-\Delta_h\phi_{\mathbf{j}}^h(t)=0,&\mathbf{j}\in\zz^d, t\in(0,T]\\\phi_{\mathbf{j}}^h(0)=\phi_{\mathbf{j}}^{h,0},\quad \partial_t\phi_{\mathbf{j}}^h(0)=\phi_{\mathbf{j}}^{h,1},&\mathbf{j}\in \zz^d.\end{array}\right.\label{DiscWave}\end{equation}
The problem (\ref{DiscWave}) is well-posed in $\dot{\hbar}^1\times\ell^2$, with
\begin{equation}\ell^2=\{\overrightarrow{f}\mbox{ s.t. }\|\overrightarrow{f}\|_{\ell^2}^2:
=h^d\sum\limits_{\mathbf{j}\in\zz^d}|f_{\mathbf{j}}|^2<\infty\} \mbox{ and } \dot{\hbar}^1=\{\overrightarrow{f}\mbox{ s.t. }\|\overrightarrow{f}\|_{\dot{\hbar}^1}^2:=\|\nabla_h^+\overrightarrow{f}\|_{\ell^2}^2=\sum\limits_{k=1}^d
\|\partial_{h,k}^+\overrightarrow{f}\|_{\ell^2}^2<\infty\}.\nonumber\end{equation}
The total energy associated to its solution is conserved in time, being defined by
$$E_h(\overrightarrow{\phi}^{h,0},\overrightarrow{\phi}^{h,1})=(\|\overrightarrow{\phi}^h(t)\|_{\dot{\hbar}^1}^2
+\|\partial_t\overrightarrow{\phi}^h(t)\|_{\ell^2}^2)/2.$$
For  a fixed $T>0$, consider the semi-discrete version of the observability inequality (\ref{ContObsIneq}):
\begin{equation}E_h(\overrightarrow{\phi}^{h,0},\overrightarrow{\phi}^{h,1})\leq C_h(T)\int_0^T\|\partial_t\overrightarrow{\phi}(t)\|_{\ell^2(\Omega^d)}^2\,dt,\mbox{ with }\|\partial_t\overrightarrow{\phi}(t)\|_{\ell^2(\Omega^d)}^2=h^d\sum\limits_{x_{\mathbf{j}}\in\Omega^d}|\partial_t\phi_{\mathbf{j}}(t)|^2.\label{DiscObsIneq}\end{equation}
For all finite $T>0$ and mesh size $h>0$, it is easy to see that (\ref{DiscObsIneq}) holds for $C_h(T)$ sufficiently large.

For the finite difference and $P_1$-classical finite element semi-discretizations of the wave equation on particular bounded domains like $d$-dimensional cubes, it is well-known from \cite{ZuaPOC} that the corresponding observability constant blows-up as $h\to 0$, because of the pathological behavior of the spurious high frequency numerical solutions. In the context of the finite difference semi-discretization of the $1-d$ wave equation on a bounded interval, the observability constant is known to blow-up exponentially as $h \to 0$ for all $T>0$ (cf. \cite{MicuBio}). As shown in \cite{LivEZDispSchr}, similar pathological high frequency phenomena have been observed in the context of the Strichartz dispersive estimates for the finite difference approximations of the Schr\"{o}dinger equation.

In view of the formulation of the observability problem under consideration, our goal here is twofold. Firstly, to show the existence of numerical waves that are concentrated in space-time and then to give its precise asymptotic form in order to illustrate the added effect of the numerical dispersion with respect to the classical continuous wave equation. These solutions propagate according to the group velocity (cf. \cite{TrefFDArt}), a notion that our construction contributes to make it precise. The group velocity of the high frequency numerical waves may be arbitrarily small, as the analysis of the dispersion diagram below shows.

Set $\Pi_h^d:=[-\pi/h,\pi/h]^d$. By applying the semi-discrete Fourier transform (SDFT) (see \cite{TrefFDArt}) on (\ref{DiscWave}) and  denoting by $\widehat{\phi}^h(\xi,t)$ the SDFT of the solution $\overrightarrow{\phi}^h(t)$ of (\ref{DiscWave}), one obtains the following  second-order ODE depending on the parameter $\xi$:
\begin{equation}\left\{\begin{array}{ll}\partial^2_t\widehat{\phi}^h(\xi,t)+\omega_{d,h}^2(\xi)\widehat{\phi}^h(\xi,t)=0,& \xi\in\Pi_h^d,t\in(0,T]\\
\widehat{\phi}^h(\xi,0)=\widehat{\phi}^{h,0}(\xi),\quad \partial_t\widehat{\phi}^h(\xi,0)=\widehat{\phi}^{h,1}(\xi),& \xi\in\Pi_h^d, \end{array}\right.\label{DiscWaveFour}\end{equation}
where $\omega_{d,h}(\xi)$ is the multi-dimensional dispersion relation associated to (\ref{DiscWave}),
$\omega_{d,h}^2(\xi)=\frac{4}{h^2}\sum_{k=1}^d\sin^2\big(\frac{\xi_k h}{2}\big)$, $\xi=(\xi_k)_{1\leq k\leq d}\in\Pi_h^d$.

The semi-discrete rays of Geometric Optics corresponding to the semi-discrete problem (\ref{DiscWave}) are straight lines of the form $x_h^{\pm}(t)=x\pm\nabla\omega_{d,h}(\xi)t$, $x\in\rr^d$, $\xi\in\Pi_h^d$, propagating with a group velocity, defined as the gradient $\nabla\omega_{d,h}(\xi)$. Note that $|\nabla\omega_{d,h}(\xi)|$ vanishes for all $\xi\in\{\pm\pi/h,0\}^d\setminus\{0\}$. This is in contrast with the behavior in the continuous case (\ref{ContWave}), where the dispersion relation is $\omega_d(\xi)=|\xi|$ and the velocity of propagation along all rays is $|\nabla\omega_d(\xi)|=1$, for all $\xi\in\rr^d$.

\section{Statement of the main result}Let $T>0$ be given. Choose an arbitrary $x^*\in B^d(0,1)$ and consider a wave number $\eta_0=h\xi_0\in\Pi_1^d\setminus\{0\}$ such that the corresponding semi-discrete ray starting at $x^*$ and traveling with velocity $|\nabla\omega_{d,h}(\xi_0)|=|\nabla\omega_{d,1}(\eta_0)|$ does not enter the observation region in time $T$, i.e. $|x^*-t\nabla\omega_{d,1}(\eta_0)|<1$, for all $t\in[0,T]$.

Consider $\gamma:=\gamma(h)>0$ such that
\begin{equation}\gamma>>1 \mbox{ and }  h\gamma<<1.\label{scale0}\end{equation}
For $\phi\in\mathcal{S}(\rr^d)$, the class of Schwartz functions in $\rr^d$,  $\widehat{\phi}$ being its continuous Fourier transform, consider the semi-discrete wave equation (\ref{DiscWave}) with initial data $(\overrightarrow{\phi}^{h,0},\overrightarrow{\phi}^{h,1})$ whose SDFT is given by \begin{equation}\widehat{\phi}^{h,0}(\xi)=\frac{1}{i\omega_{d,h}(\xi)}\sqrt{\frac{2\pi}{\gamma}}^d\widehat{\phi}\left(\frac{\xi-\xi_0}{\gamma}\right)
\exp(-ix^*\cdot(\xi-\xi_0))\chi_{\Pi_h^d}(\xi)\mbox{ and }\ \widehat{\phi}^{h,1}(\xi)=i\omega_{d,h}(\xi)\widehat{\phi}^{h,0}(\xi).\label{initialdataSylvain}\end{equation}

These data, in the physical space, correspond, roughly, to highly oscillating Gaussian profiles. As shown by S. Ervedoza in \cite{Erv}, by a stationary phase like argument, one can show that for all $\beta\in\rr_{+}$, there exists a constant $C_{\beta}(\widehat{\phi},T)$ such that the observability constant $C_h(T)$ in (\ref{DiscObsIneq}) satisfies $C_h(T)\geqslant C_{\beta}(\widehat{\phi},T)\gamma^{\beta}$.

The initial data (\ref{initialdataSylvain}) are concentrated high frequency wave packets at the wave number $\xi=\xi_0=\eta_0/h$ with a width of order $1/\gamma$, with $\gamma<<1/h\to \infty$ as $h\to 0$. This limits drastically the spread of the wave packet and allows to concentrate its energy around the chosen ray and to reduce the added dispersive effects that the non-trivial Hessian matrix of the dispersion relation introduces.  Simultaneously, by the uncertainty principle, this forces the wave packet, in the physical space, to have a spread factor of order  $1/\gamma>>h$, which is asymptotically larger than the characteristic size $h$ of the mesh. This is natural, in the sense that the numerical effects are only detected when an infinite number of nodes enter asymptotically in the determination of the data of the solution.

This example of wave packet clearly illustrates the classical effect due to the group velocity that is plotted in Figure \ref{Fig2} below. We plot in red the initial velocity $\phi^1(x)=\exp(-\gamma|x|^2/2)\exp(i\xi_0 x)$ (we also take $i|\xi|\widehat{\phi^0}=\widehat{\phi^1}$), with $h=0.005$, $\xi_0 h=19\pi/20$, $\gamma=h^{-\alpha}$, $2\alpha=0.75$ in $d=1$ (left) and $d=2$ (right), which coincides, up to an exponentially small error (with respect to $h$), with the discrete initial data (\ref{initialdataSylvain}). In blue, we plot the time derivative of the solution of the continuous wave equation (\ref{ContWave}) and in black the one corresponding to the semi-discrete wave equation (\ref{DiscWave}) with initial data (\ref{initialdataSylvain}) at time $t=1$. For $d=2$, we represent only the projection of these (continuous and discrete) solutions on the $x$-plane. This experiment shows that, as theory predicts, the semi-discrete wave packets propagate with a velocity that is much smaller than the one of the continuous wave equation.
\begin{figure}[h!]
\begin{center}
  \includegraphics[width=6cm,height=4cm]{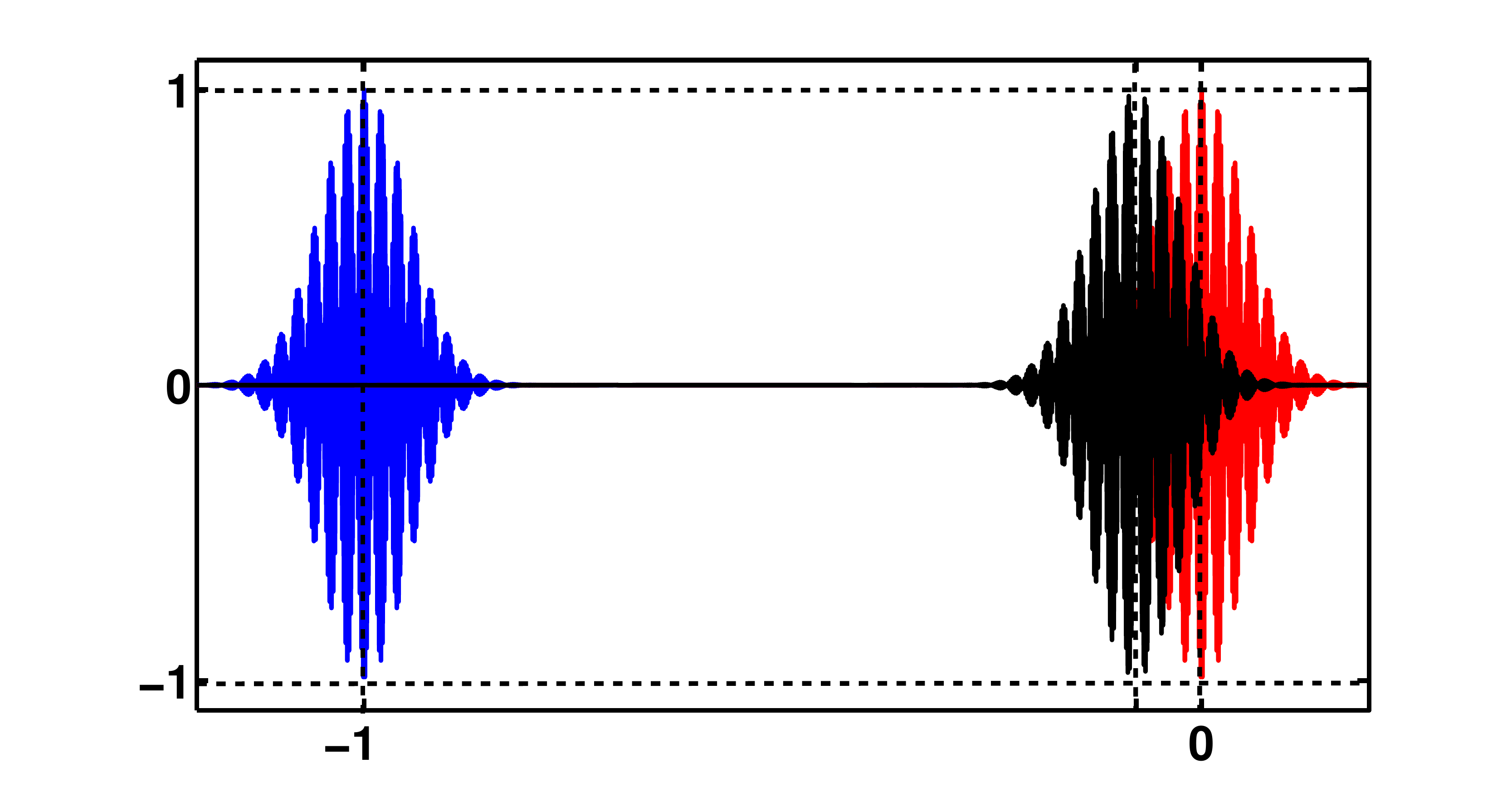}\includegraphics[width=6cm,height=4cm]{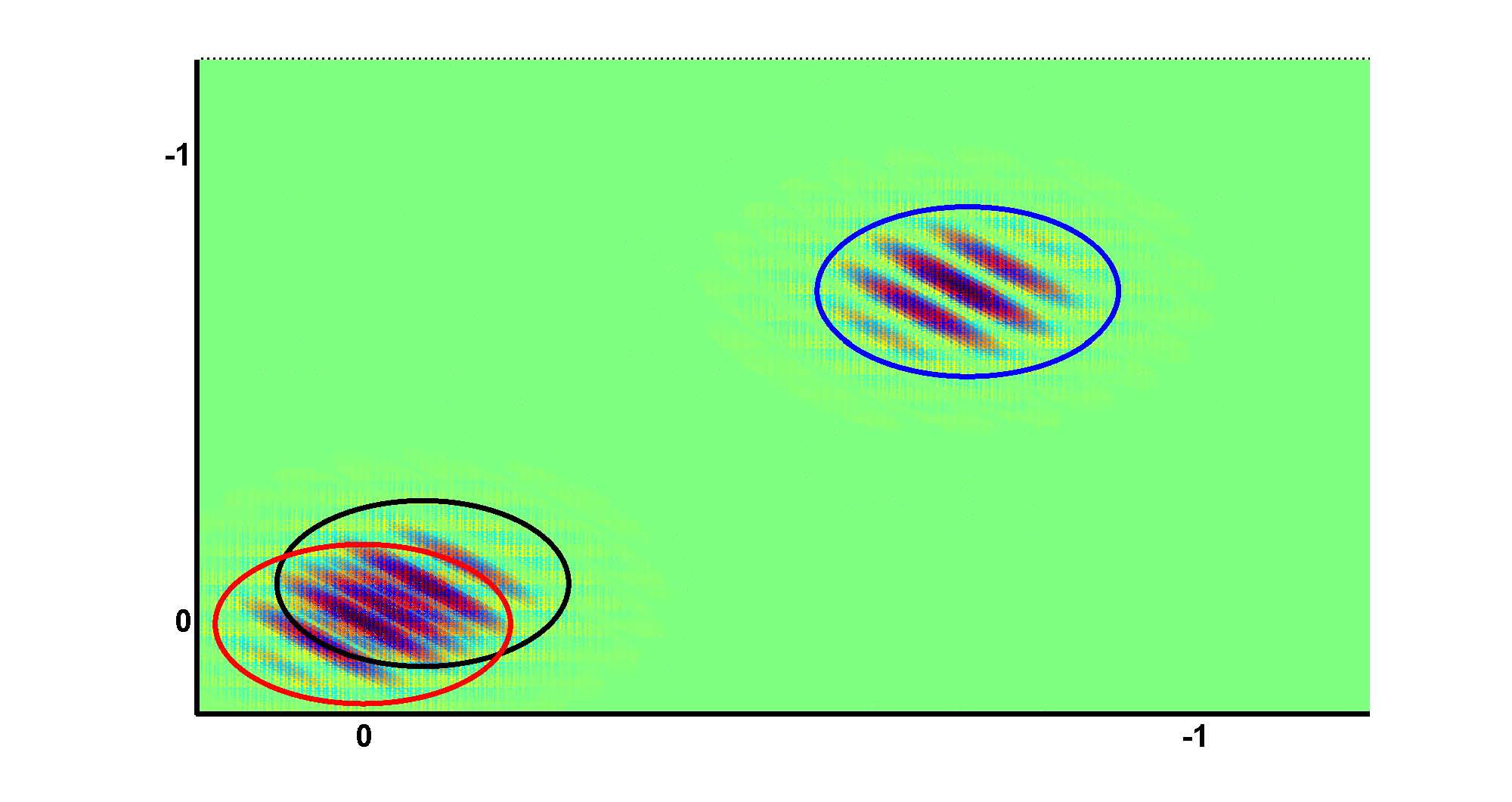}\\
  \caption{Transport of continuous versus discrete high frequency wave packets in dimension $d=1$ and $d=2$.}
  \label{Fig2}
  \end{center}
\end{figure}

In the following, we analyze more precisely the behavior of these wave packets. To do this, we appro\-xi\-mate the analytic dispersion relation $\omega_{d,h}(\xi)$ by a polynomial one.

Set $\Pi_h^d-\xi_0:=\{\xi-\xi_0,\xi\in\Pi_h^d\}$. The time derivative of the solution of (\ref{DiscWave}) with initial data (\ref{initialdataSylvain}) is given by the following wave packet
\begin{equation}\partial_t\phi_{\mathbf{j}}(t)=\frac{1}{(2\pi)^d}\int\limits_{\rr^d}\sqrt{\frac{2\pi}{\gamma}}^d\widehat{\phi}(\xi/\gamma)
\chi_{\Pi_h^d-\xi_0}(\xi)\exp(it\omega_{d,h}(\xi+\xi_0))
\exp(i\xi\cdot(x_{\mathbf{j}}-x^*))\,d\xi.\label{solnnumerica}\end{equation}

To simplify the presentation and, without loss of generality, we set $\xi_0=0$ and $x^*=0$. Since $\widehat{\phi}\in \mathcal{S}(\rr^d)$, we may also neglect the characteristic function in (\ref{solnnumerica}). Setting $\omega(\xi)=\omega_{d,1}(\xi+\eta_0)$ and $x=x_{\mathbf{j}}-x^*$, the wave packet (\ref{solnnumerica}) can be written as
\begin{equation}u(x,t)=\frac{1}{(2\pi)^d}\int\limits_{\rr^d}\sqrt{\frac{2\pi}{\gamma}}^d\widehat{\phi}\big(\frac{\xi}{\gamma}\big)
\exp\Big(\frac{it}{h}\omega(\xi h)\Big)\exp(i\xi\cdot x)\,d\xi.\label{solnnumerica2}\end{equation}

In view of the analyticity of $\omega(\eta)$, we can split it as $\omega(\eta)=L(\eta)+D(\eta)+R(\eta), $ where $L(\eta)=\omega(0)+\nabla\omega(0)\cdot\eta$ is the linear part of $\omega$, $D$ is the second order term in the Taylor expansion of $\omega$ about $\eta=0$ and $R$ is the corresponding reminder, given explicitly by
$$D(\eta)=\sum\limits_{|\alpha|=2}\frac{1}{\alpha!}D^{\alpha}\omega(0)\eta^{\alpha},
R(\eta)=\sum\limits_{|\alpha|=3}\frac{3}{\alpha!}\eta^{\alpha}\int\limits_0^1(1-\lambda)^2D^{\alpha}\omega(\lambda\eta)\,d\lambda.$$

Factoring out the time dependent complex exponential generated by the zero order term in the Taylor expansion of the dispersion relation, we may decompose $u$ as $u\exp(-it\omega_{d,1}(0)/h)=v+v^{R}$, where
$$v(x,t)=\frac{1}{(2\pi)^d}\int\limits_{\rr^d}\sqrt{\frac{2\pi}{\gamma}}^d\widehat{\phi}\big(\frac{\xi}{\gamma}\big)\exp\Big(\frac{it}{h}D(\xi h)\Big)\exp(i\xi\cdot (x+t\nabla\omega(0)))\,d\xi$$
and
$$v^R(x,t)=\frac{1}{(2\pi)^d}\int\limits_{\rr^d}\sqrt{\frac{2\pi}{\gamma}}^d\widehat{\phi}\big(\frac{\xi}{\gamma}\big)\exp\Big(\frac{it}{h}D(\xi h)\Big)\Big(\exp\Big(\frac{it}{h}R(\xi h)\Big)-1\Big)\exp(i\xi\cdot (x+t\nabla\omega(0)))\,d\xi.$$

Our main result determines the conditions on $\gamma$ and $\widehat{\phi}$, such that $v^R$ is a reminder term.
\begin{theorem}For all $t>0$ and $\widehat{\phi}\in \mathcal{S}(\rr^d)$, the following estimate holds:
\begin{equation}\frac{\|v^{R}(\cdot,t)\|_{L^2(\rr^d)}^2}{\|u(\cdot,t)\|_{L^2(\rr^d)}^2}\leq C(\widehat{\phi})h^{4}\gamma^{6}t^2,\label{estimRk}\end{equation}
where $C(\widehat{\phi})=C\||\cdot|^3\widehat{\phi}\|_{L^2(\rr^d)}^2/\|\widehat{\phi}\|_{L^2(\rr^d)}^2$ and $C=\Big(\sum\limits_{|\alpha|=3}\frac{1}{\alpha!}\|D^{\alpha}\omega\|_{L^{\infty}(B(0,h\gamma))}\Big)^2$
\end{theorem}
For a finite time interval $t\in[0,T]$, one can guarantee that $v^{R}$ is small in the sense of (\ref{estimRk}) when
\begin{equation}\gamma h^{2/3}<<1,\label{scale1}\end{equation} which is a more restrictive condition on $\gamma$ than (\ref{scale0}), guaranteing that the energy concentrated outside the ray is polynomially small. Indeed, (\ref{estimRk}) yields an asymptotic description of the solution globally in space-time, i.e. $u\thicksim\exp(it\omega(0)/h)v$.
Intuitively, the scale (\ref{scale1}) is motivated by the fact that $|R(\xi h) /h | \le \sqrt{C} h^2 |\xi|^3$. For this to be asymptotically small, the width of the Fourier transform of the profile has to be limited by (\ref{scale1}).

The function $v$ is a solution of the following PDE:
\begin{equation}\partial_tv=\nabla\omega(0)\cdot\nabla_xv-ih\sum\limits_{|\alpha|=2}D^{\alpha}\omega(0)D^{\alpha}_xv.\label{eqn1}\end{equation}
This is a transport equation perturbed by an asymptotically small (as $h$ tends to zero) Schr\"{o}dinger like second order term. The following result emphasizes that there exists solutions of (\ref{eqn1}) for which the relevant scale is $$\gamma h^{1/2}=1.$$
\begin{theorem}There exist solutions of (\ref{eqn1}) admitting the following asymptotic expansion:
\begin{equation}v(x,t)=\sum\limits_{j=0}^{\infty}h^{j/2}a_j\left(\frac{x}{h^{1/2}},\frac{t}{h^{1/2}}\right)
\exp\left(i\xi\cdot\frac{\eta_0}{h^{1/2}}+it\omega_2\left(\frac{\eta_0}{h^{1/2}}\right)\right),\label{expansion}\end{equation}
where $\eta_0\in\rr^d$ is a fixed wave number, $\omega_2(\eta)=\nabla\omega(0)\cdot\eta+h\sum_{|\alpha|=2}\frac{1}{\alpha!}D^{\alpha}\omega(0)\eta^{\alpha}$ is the dispersion relation corresponding to (\ref{eqn1}) and $(a_j(x,t))_{j\in\nn}$ solve the following system of PDEs:
$$\partial_ta_0=\nabla\omega(0)\cdot\nabla_xa_0\mbox{ and }\partial_ta_{j+1}=\nabla\omega(0)\cdot\nabla_xa_{j+1}-i\sum\limits_{|\alpha|=2}\frac{1}{\alpha!}D^{\alpha}\omega(0)\sum\limits_{0\leq\beta<\alpha}
(i\eta_0)^{\beta}D^{\alpha-\beta}_xa_j,\forall j\in\nn.$$
\label{teoasympk3}\end{theorem}
By taking the support of the initial datum of $a_0$ to be compact, one can observe that $v$ is concentrated along a neighborhood of the ray of width $\sqrt{h}$. This scale $\sqrt{h}$ is critical due to the added dispersion that the Schr\"{o}dinger like term introduces.  Once the asymptotic expansion of $v$ is given as in (\ref{expansion}), one immediately gets that of $u$ and therefore of $\partial_t\overrightarrow{\phi}$ and of $\overrightarrow{\phi}$. In this way, we get the asymptotic form of the high frequency wave packet.

This kind of expansion can be further developed, incorporating higher order terms of the Taylor expansion of the dispersion relation and a multiple-scale ansatz. This issue will be developed in detail in \cite{ErvMarZua}.

\section*{Acknowledgements}
The authors acknowledge Sylvain Ervedoza for fruitful discussions. Both authors were partially supported by the Grant MTM2008-03541 of the MICINN, Spain, and the ERC Advanced Grant FP7-246775 NUMERIWAVES.

\end{document}